\documentclass[a4paper,10pt]{article}
\usepackage[pdftex]{graphicx}
\DeclareGraphicsExtensions{.eps}

\usepackage{newtxtext,newtxmath} 
\usepackage{amsmath,cite,subcaption,multirow,color,url,algorithm}
\usepackage[noend]{algpseudocode}
\usepackage[vmargin={2cm,2cm},hmargin=3cm]{geometry}

\title{Routing optimization on power packet dispatching system based on energy loss minimization}
\author{Shiu Mochiyama\thanks{Email: s-mochiyama@dove.kuee.kyoto-u.ac.jp} \and Kazuhiro Koto \and Takashi Hikihara}
\date{Department of Electrical Engineering, \\Kyoto University}

\begin{document}
\maketitle

\begin{abstract}
 Power packet dispatching system has been proposed for smart power management in the form of discretized packet. 
 In this paper, we discuss the routing optimization of power packets on the network of power routers. 
 We propose a cost metric for the power packet delivery by circuit analysis of the router network. 
 Using the metric, we formulate the optimization problem as a general shortest path problem from a source node to a load node. 
 The result of numerical simulations shows that the proposed algorithm can allocate distributed power sources to load demands and identify the optimal path for the power delivery. 
\end{abstract}

\section{Introduction}
Electric energy management has been an essential technology in a variety of systems such as vehicles, robots, and microgrids \cite{Chen.etal-2018,Seok.etal-2015,Yang.etal-2018}. 
Particularly in such systems, the strong demand for eco-friendliness has raised the importance of an efficient inclusion of distributed power sources. 
Then, the time-varying profile of the sources creates more complex power flow than in the conventional system based on a large stable power source. 
A key enabler for handling the complex power flow is the support of Information and Communication Technology (ICT) \cite{Kim.Kumar-2012,He.etal-2008,Gelenbe.Abdelrahman-2018}.

Power packet dispatching system has been proposed as a way of such smart power management \cite{Takuno.etal-2010,Takahashi.etal-2015}. 
\begin{figure}[tb]
 \centering
 \includegraphics[width=.5\columnwidth]{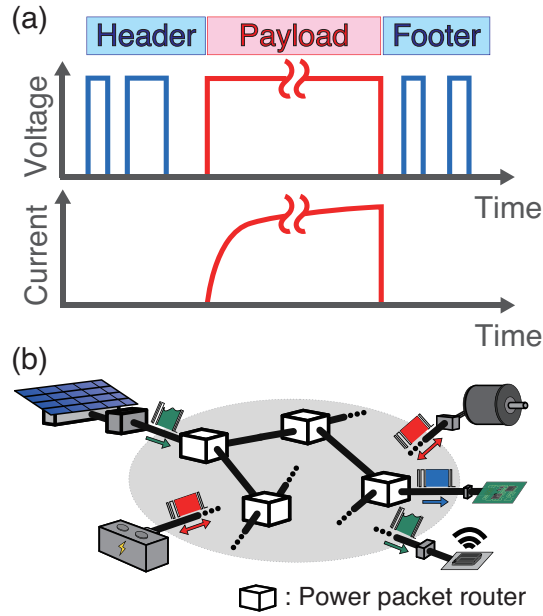}
 \caption{Power packet dispatching system\cite{Takuno.etal-2010,Mochiyama.etal-2021}. }
 \label{fig:ppds}
\end{figure}
Figure~\ref{fig:ppds} depicts the concept of the power packet dispatching system. 
A power packet is a unit of transfer of both power and information. 
A payload delivers a quantized unit of electric power. 
An information tag, attached to the payload in voltage waveform, represents information such as identification of the packet. 
The physical tag attachment enables the network to distinguish each packet and thus to process electric power in a fully digitized manner. 
This feature makes the power management highly compatible with the ICT.

A power packet is transferred from its source to a load through a network of apparatuses called power routers (hereafter we simply call them routers). 
Equipped with some temporary energy storage, the routers deliver a power packet in a store and forward manner \cite{Takahashi.etal-2015}. 
Recent study \cite{Mochiyama.etal-2021} reported that the power packetization can be realized at a rating of several kW.

A remaining challenge is a routing optimization of a power packet. 
Similarly with the IP routing \cite{Fortz.etal-2002}, in a power packet dispatching system of a number of routers, there are many possible paths between a pair of a source to a load. 
It then becomes necessary to set an appropriate strategy for selecting one. 
However, we cannot directly apply the routing problem for the IP routing. 
As we consider the flow of physical quantity, we need a new criterion for evaluating the power packet routing based on the dynamical behavior of the circuit.

In addition, a matching optimization of the distributed sources and loads is another important challenge. 
The network structure of a power packet dispatching system can vary dynamically. 
In a microgrid application, for example, electric vehicles (EVs) can connect to or disconnect from the system at an arbitrary time. 
Renewable sources such as photovoltaic cells and wind turbines can also have time-dependent output profiles. 
In this way, the sources are spatio-temporally distributed. 
Now suppose that a load does not have a specific demand for the \textit{origin} of power but only for the \textit{amount} of power, the system is requested to select one of the multiple sources available at the time.

In this paper, we address the above challenges by introducing an optimization of power packet routing based on a cost metric derived by circuit analysis. 
First, we propose a circuit model of power packet transfer.
The model lets us derive the cost of power packet transfer between arbitrary node pairs on the dispatching network. 
Second, with the derived cost metric, we consider the path optimization problem. 
We develop an algorithm to solve the problem by reducing it to a general shortest path problem on a graph. 
Lastly, we verify the proposed algorithm through numerical simulations.

\section{Model of power packet transfer}
\begin{figure}[tb]
 \centering
 \includegraphics[width=.5\columnwidth]{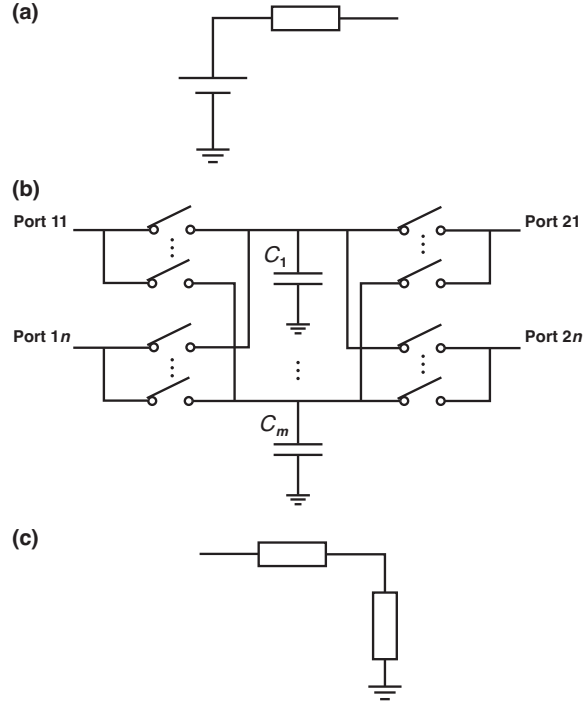}
 \caption{Circuit models for (a) source, (b) router, and (c) load.}
 \label{fig:elem}
\end{figure}

We assume a power packet dispatching system consisting of $n_\mathrm{s}$ sources, $n_\mathrm{r}$ routers, and $n_\mathrm{l}$ loads, where $n_\mathrm{s}$, $n_\mathrm{r}$, and $n_\mathrm{l}$ are positive integer. 
Figure~\ref{fig:elem} shows the circuit models for a source, a router, and a load. 
A source consists of an ideal voltage source, an ideal switch, and series resistance. 
Turning on and off the switch, a source produces a power packet. 
The series resistance represents whole resistance between the voltage source and the output port, including internal resistance, switch's conduction resistance, and line resistance. 
A router is equipped with multiple input/output ports that accept bidirectional power flow. 
The circuit consists of ideal switches, capacitors, and series resistance defined similarly to a source.
We denote the number of capacitor in a router by $n_\mathrm{c}$. 
The capacitors are used as a temporary power storage in receiving a power packet and as a source in forwarding a power packet.

Now the dispatching system can be represented by a simple graph $G=G(\mathcal{N},\mathcal{E})$, where $\mathcal{N}$ and $\mathcal{E}$ represent sets of nodes (sources, routers, and loads) and edges (power lines between nodes), respectively. 
For each node $i \in \mathcal{N}$ we set an associated attribute $a_i$ that represents the role of the node, i.e. $a_i \in \{\mathrm{s},\mathrm{r},\mathrm{l}\}$, where the elements in the set represent ``source,'' ``router,'' and ``load,'' respectively. 
For each edge $e_{(i,j)} = (i,j) \in \mathcal{E}$, where $i,j \in \mathcal{N}$, we set an associated attribute $w_{(i,j)} (> 0)$ that represents a \textit{cost} for power delivery from $i$ to $j$. 
How to determine the cost is arbitrary at this stage;
we will define one candidate later in this paper, with which we conduct numerical simulations.

Here we also define a state variable of the graph, namely voltage of each node: 
\begin{equation}
 \boldsymbol{v} = 
  \begin{bmatrix}
   \boldsymbol{v}_1 \\
   \vdots \\
   \boldsymbol{v}_N
  \end{bmatrix}
\end{equation}
where $N := n_\mathrm{s} + n_\mathrm{r} + n_\mathrm{l}$ and the row vectors $\boldsymbol{v}_i \in \mathbb{R}^{m_\mathrm{c}}$ ($i=1,2,\dots,N$) represents voltage of node $i \in N$ in the following manner. 
\begin{itemize}
 \item If $a_i = \mathrm{r}$, $v_i = [v_{(i,1)},\cdots,v_{(i,m_\mathrm{c})}]$ (voltage of the capacitors of the node).
 \item If $a_i = \mathrm{s}$, $v_i = [V_i,\cdots,V_i]$, where $V_i$ is voltage of the ideal source (filled with the voltage of the source).
 \item If $a_i = \mathrm{l}$, $v_i = [0,\cdots,0]$ (filled with zero).
\end{itemize}
Note that the duplication of elements when $a_i = \mathrm{s}$ or $\mathrm{l}$ is just for keeping the size of the matrix.

Here we set some assumptions on power delivery within the dispatching system. 
The first one is on the time synchronization of the system. 
We fix the time duration of one power packet%
\footnote{Note that we do \textit{not} fix the time duration of payload.}
at $T$, and denote the $k$-th time slot, $(k-1)T \leq t < kT$, simply by $k$ ($k=1,2,\dots$), where $t$ represents continuous time. 
Second, we fix the amount of energy transferred in a time slot on each edge and denote it by $U$. 
This is achieved by changing the length of the payload. 
This restriction is to minimize the energy change in the network. 
Third, when a router is connected to another node, no overlap of capacitor is allowed. 
For example, suppose that a router node $i$ receives a power packet from $j$ and forwards it to $k$, the capacitors used for delivery on $(j,i)$ and on $(i,k)$ must be different. 
This restriction is to ensure a one-to-one delivery of unit energy.

\begin{figure}[tb]
 \centering
 \includegraphics[width=.33\columnwidth]{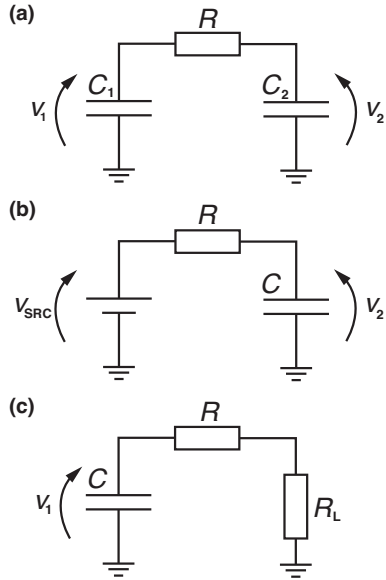}
 \caption{Circuit models for the three types of connection; (a) Router-router connection, (b) Source-router connection, and (c) Router-load connection. }
 \label{fig:circuit}
\end{figure}

Under the assumptions, the connection between any nodes can be classified into three types: router-router, source-router, and router-load connections. 
Figure~\ref{fig:circuit} shows the circuit model for the three types of connection. 
The detailed operation of each connection is explained in the following subsections.

\subsection{Router-router Connection}
In this connection, one of the capacitors in the router forwards power packet to one of the capacitors in the other router (Fig.~\ref{fig:circuit}~(a)). 
For simplicity, resistance of all components between the capacitors including the switches and the lines is represented by $R$. 

Now, letting $v_1(k)$ and $v_2(k)$ be voltage of left and right side capacitors at discrete time $k$, respectively, we derive their update law for one power packet transfer, i.e. $v_1(k+1)$ and $v_2(k+1)$ as functions of $v_1(k)$ and $v_2(k)$. 
We assume $C_1 = C_2$ in the rest of this paper. 
This is just for visibility of derived equations and results of numerical simulation conducted later; they can, of course, have different values without any problem in derivation or simulation. 
Applying Kirchhoff's circuit laws and carrying out the time integration, we have
\begin{align}
 v_1(k+1) &= \frac{1}{2} \left\{ (v_1(k)-v_2(k)) \mathrm{e}^{\frac{-2T_0}{CR}} + v_1(k) + v_2(k) \right\} , \label{eq:upd_rr1} \\
 v_2(k+1) &= \frac{1}{2} \left\{-(v_1(k)-v_2(k)) \mathrm{e}^{\frac{-2T_0}{CR}} + v_1(k) + v_2(k) \right\} , \label{eq:upd_rr2} 
\end{align}
where $T_0$ ($\leq T$) represents the payload duration. 

Here is derived the amount of energy exchange in a power packet as follows: 
\begin{align} 
 \Delta E_\mathrm{send} &= \frac{1}{8}C (v_1(k)-v_2(k)) (1- \mathrm{e}^{\frac{-2T_0}{CR}}) \left\{ (3v_1(k)+v_2(k)) + (v_1(k)-v_2(k)) \mathrm{e}^{\frac{-2T_0}{CR}} \right\} , \\
 \Delta E_\mathrm{receive} &= \frac{1}{8}C (v_1(k)-v_2(k)) (1- \mathrm{e}^{\frac{-2T_0}{CR}}) \left\{ (v_1(k)+3v_2(k)) - (v_1(k)-v_2(k)) \mathrm{e}^{\frac{-2T_0}{CR}} \right\} , \label{eq:rr_Er} \\
 \Delta E_\mathrm{loss} &= \frac{1}{4}C (v_1(k)-v_2(k))^2 \left(1 - \mathrm{e}^{\frac{-4T_0}{CR}}\right) , \label{eq:rr_El}
\end{align}
where $\Delta E_\mathrm{send}$, $\Delta E_\mathrm{receive}$, and $\Delta E_\mathrm{loss}$ represent the amount of energy sent from the left side router, received in the right side router, and lost as Joule heat, respectively. 

Using Eq.~(\ref{eq:rr_Er}), we can estimate the packet length required to deliver unit energy $U$. 
For $T \ll CR$, the linearization of Eq.~(\ref{eq:rr_Er}) gives
\begin{equation}
 \tau = \frac{RU}{ \{v_1(k)-v_2(k) \} v_2(k)} ,
  \label{eq:pleng_rr}
\end{equation}
where $\tau$ denotes the estimated value of $T_0$.

\subsection{Source-router Connection}
In this connection, the voltage source supplies one of the capacitors in the receiving router (Fig.~\ref{fig:circuit}~(b)). 
Similarly with the router-router connection, $R$ represents all resistive components between the nodes. 

Letting $v_\mathrm{src}$ and $v_2(k)$ be voltage of the source and the right side capacitor at discrete time $k$, respectively, we have
\begin{equation}
 v_2(k+1) = (v_2(k) - v_\mathrm{src}) \mathrm{e}^{\frac{-T_0}{CR}} + v_\mathrm{src}. 
  \label{eq:upd_sr} 
\end{equation}

The amount of energy transfer can be derived as follows:
\begin{align}
 \Delta E_\mathrm{send} &= C v_\mathrm{src} (v_\mathrm{src} - v_2(k)) \left(1 - \mathrm{e}^{\frac{-T_0}{CR}} \right) , \\
 \Delta E_\mathrm{receive} &= \frac{1}{2} C (v_\mathrm{src} - v_2(k)) \left(1 - \mathrm{e}^{\frac{-T_0}{CR}} \right) \left\{ v_\mathrm{src} \left(1 - \mathrm{e}^{\frac{-T_0}{CR}} \right) + v_2(k) \left(1 + \mathrm{e}^{\frac{-T_0}{CR}} \right)\right\} , \label{eq:sr_Er} \\
 \Delta E_\mathrm{loss} &= \frac{1}{2} C (v_\mathrm{src} - v_2(k))^2 \left(1 - \mathrm{e}^{\frac{-2T_0}{CR}} \right) . \label{eq:sr_El}
\end{align}

For $T \ll CR$, the linearization of Eq.~(\ref{eq:sr_Er}) gives
\begin{equation}
 \tau = \frac{RU}{ \{ v_\mathrm{src}-v_2(k) \} v_2(k)}.
  \label{eq:pleng_sr}
\end{equation}

\subsection{Router-load Connection}
In this connection, one of the capacitors in the sender router supplies the load (Fig.~\ref{fig:circuit}~(c)). 
Here $R$ represents the resistance of all components but the load. 
The load resistance is denoted by $R_\mathrm{L}$. 

Letting $v_1(k)$ be voltage of the capacitor at discrete time $k$, we have
\begin{equation}
 v_1(k+1) = v_1(k) \mathrm{e}^{\frac{-T}{C(R+R_\mathrm{L})}} .
\label{eq:upd_rl} 
\end{equation}

The amount of energy transfer can be derived as follows:
\begin{align}
 \Delta E_\mathrm{send} &= \frac{1}{2} C v^2_{01} \left( 1 - \mathrm{e}^{\frac{-2T}{C(R+R_\mathrm{L})}} \right) , \\
 \Delta E_\mathrm{receive} &= \frac{1}{2} \frac{R_\mathrm{L}}{R + R_\mathrm{L}} C v^2_{01} \left( 1 - \mathrm{e}^{\frac{-2T}{C(R+R_\mathrm{L})}} \right) , \label{eq:rl_Er} \\
 \Delta E_\mathrm{loss} &= \frac{1}{2} \frac{R}{R + R_\mathrm{L}} C v^2_{01} \left( 1 - \mathrm{e}^{\frac{-2T}{C(R+R_\mathrm{L})}} \right) . \label{eq:rl_El}
\end{align}

For $T \ll C(R+R_\mathrm{L})$, the linearization of Eq.~(\ref{eq:rl_Er}) gives
\begin{equation}
 \tau = \frac{(R+R_\mathrm{L})^2 U}{R_\mathrm{L}v^2_1(k)} .
  \label{eq:pleng_rl}
\end{equation}

\section{Flow Control}
Based on the circuit models derived above, we consider flow control on the power packet dispatching system.
A general router network contains multiple distributed sources available. 
Thus, considering a power supply to a specific load, there are multiple candidates of the power source to be used. 
In this section, we develop an algorithm to find the optimal choice of the source and the path from the source to the load on the router network. 

Apparently there are arbitrary definition for being the \textit{optimal}. 
In this paper, we define it as ``the energy loss becomes the lowest among the possible paths.'' 
This is a natural way of defining the cost because energy loss is one of the most critical criteria in a power distribution network. 

We define the cost by the proportion of the energy loss to the total energy received, namely $\Delta E_\mathrm{loss}/\Delta E_\mathrm{receive}$. 
With Eqs.~(\ref{eq:rr_Er})~and~(\ref{eq:rr_El}), we can calculate the cost for router-router connection as follows:
\begin{equation}
 w = \frac{-2(v_0(k)-1)(1- \mathrm{e}^{\frac{-4T}{CR}})}{(v_0(k)-1) \mathrm{e}^{\frac{-4T}{CR}} + 2(v_0(k)+1)\mathrm{e}^{\frac{-2T}{CR}} - 3_0(k)v - 1 }, 
\end{equation}
where $v_0(k) = v_1(k)/v_2(k)$. The definitions of $v_1$, $v_2$, $T$, $C$, and $R$ are same as introduced in Eqs.~(\ref{eq:rr_Er})~and~(\ref{eq:rr_El}). 
Similarly, for source-router connection, we have
\begin{equation}
 w = \frac{v-1}{2v} \left( 1 + \mathrm{e}^{\frac{-T}{CR}} \right) ,
\end{equation}
and for router-load connection, we have
\begin{equation}
 w = \frac{R}{R + R_\mathrm{L}} . 
\end{equation}

Here we define the path.  
Given a pair of a source node $n_{\mathrm{s}} \in \mathcal{N}$ and a load node $n_{\mathrm{d}} \in \mathcal{N}$ of a graph $G = G(\mathcal{N},\mathcal{E})$, a path from $n_{\mathrm{s}}$ to $n_{\mathrm{d}}$, $p(n_{\mathrm{s}},n_{\mathrm{d}})$, is a sequence of edges that connects the source and load without loop.
That is, $p(n_{\mathrm{s}},n_{\mathrm{d}}) = ((n_1,n_2),(n_2,n_3),\dots,(n_{k-1},n_k))$ s.t. $n_i \in \mathcal{N}$ ($i = 1,2,\dots,k$), $(n_j,n_{j+1}) \in \mathcal{E}$ ($j=1,2,\dots,k-1$), $n_1 = n_{\mathrm{s}}$, $n_k = n_{\mathrm{d}}$, and $n_i \neq n_j$ ($i \neq j$).

Based on the definitions above, the problem to be solved is expressed as follows.
Given a source node $n_{\mathrm{s}}$ and a load node $n_{\mathrm{d}}$, find $p(n_{\mathrm{s}},n_{\mathrm{d}})$ along which the sum of energy loss $\sum_{i} w_{(n_i,n_{i+1})}$ becomes the lowest.

\begin{algorithm}[tbp]
  \caption{Pseudo code for flow control. }\label{alg:flow}
   \begin{algorithmic}[1]
    \Procedure{CalcOnePacketTransfer}{$G,\boldsymbol{v}$}
    \For{each edge $(i,j)$ in $E$}
     \If{$a_i == r$ or $a_j == r$}
      \State $(c_i,c_j) \gets \textsc{selectCap}()$
      \Comment $(c_i,c_j)$: indexes of capacitors used on $(i,j)$
     \EndIf
     \If{$(a_i,a_j) == (r,r)$ and $v_{(i,c_i)} > v_{(j,c_j)}$}
      \State $w_{(i,j)} \gets \textsc{calcEnergyLossRR}(v_{(i,c_i)},v_{(j,c_j)})$
      \Comment Apply Eq.~(\ref{eq:rr_El})
     \ElsIf{$(a_i,a_j) == (s,r)$ and $v_i > v_{(j,c_j)}$}
      \State $w_{(i,j)} \gets \textsc{calcEnergyLossSR}(V_i,v_{(j,c_j)})$
      \Comment Apply Eq.~(\ref{eq:sr_El})
     \ElsIf{$(a_i,a_j) == (r,l)$}
      \State $w_{(i,j)} \gets \textsc{calcEnergyLossRL}(v_{(i,c_i)},0)$
      \Comment Apply Eq.~(\ref{eq:rl_El})
     \Else
      \State $w_{(i,j)} \gets \infty$
     \EndIf
    \EndFor
    \Statex
    \State $p \gets \textsc{solveSPF}(W,n_\mathrm{s},n_\mathrm{d})$
    \Comment Apply shortest path solver
    \Statex
    \For{each edge $(i,j)$ in $p$}
     \If{$(a_i,a_j) == (r,r)$}
      \State $\tau_{(i,j)} \gets \textsc{calcPayloadLengRR}(i,j,v_{(i,c_i)},v_{(j,c_j)};U)$
      \Comment Apply Eq.~(\ref{eq:pleng_rr})
     \ElsIf{$(a_i,a_j) == (s,r)$}
      \State $\tau_{(i,j)} \gets \textsc{calcPayloadLengSR}(i,j,V_i,v_{(j,c_j)};U)$
      \Comment Apply Eq.~(\ref{eq:pleng_sr})
     \ElsIf{$(a_i,a_j) == (r,l)$}
      \State $\tau_{(i,j)} \gets \textsc{calcPayloadLengRL}(i,j,v_{(i,c_i)},0;U)$
      \Comment Apply Eq.~(\ref{eq:pleng_rl})
     \EndIf
    \EndFor
    \Statex
    \For{each edge $(i,j)$ in $p$}
     \If{$(a_i,a_j) == (r,r)$}
      \State $v_{(i,c_i)}, v_{(j,c_j)} \gets \textsc{updateVoltRR}(\tau_{(i,j)},v_{(i,c_i)},v_{(j,c_j)})$
      \Comment Apply Eqs.~(\ref{eq:upd_rr1})~and~(\ref{eq:upd_rr2})
     \ElsIf{$(a_i,a_j) == (s,r)$}
      \State $v_{(j,c_j)} \gets \textsc{updateVoltSR}(\tau_{(i,j)},V_i,v_{(j,c_j)})$
      \Comment Apply Eq.~(\ref{eq:upd_sr})
     \ElsIf{$(a_i,a_j) == (r,l)$}
      \State $v_{(i,c_i)} \gets \textsc{updateVoltRL}(\tau_{(i,j)},v_{(i,c_i)},0)$
      \Comment Apply Eq.~(\ref{eq:upd_rl})
     \EndIf
    \EndFor
    \EndProcedure
   \end{algorithmic}
\end{algorithm}
The algorithm including path finding and update of voltage variables is shown in Algorithm~{\ref{alg:flow}}.
As the first step, it updates the cost between all possible connections of nodes (Lines 2--12).
At the beginning of this step, we choose one capacitor of those in each router to avoid duplicated use (Line 2).
Then, the cost on all the edges is calculated. 
Next, the second step is the shortest path finding (Line 13). 
This can be achieved by any algorithm for all pair shortest path problem.
We apply Floyd–Warshall algorithm \cite{Gallo.Pallottino-1988} in the following numerical simulation. 
The third step is calculation of payload length $T_0$ (Lines 14--15). 
This step ensures the unit energy transfer along the path. 
The last step is to update voltage of capacitors used in the path (Lines 16--22).

Here is given a comment on the step of capacitor selection (Line 2). 
The selection of capacitors affects the result of the path finding.
In other words, the algorithm gives the optimal path under a particular selection of the capacitor indexes. 
In this study we avoid its arbitrariness by setting the initial voltage of all capacitors within a router at a uniform value. 
If they are at different voltage levels, one way is to select the pair with lowest ratio of voltage difference.

\section{Numerical Simulation}

\subsection{Methods}
We used python 3.9.1 and the module ``shortest\_path'' of SciPy 1.6.2 \cite{TheSciPycommunity-} for implementing the algorithm.

We set the following simplifications in the setup in order to put the focus of the verification on the essential targets of this paper (see Introduction). 
\begin{itemize}
 \item The same values were set for $C$, $R$, and $U$, $C = 1.0\times 10^{-} \,\mathrm{F}$, $R = 1.0\times 10^{-} \,\mathrm{\Omega}$, and $U_0 = 3.0 \times 10^{-3} \,\mathrm{J}$, for all the connections. 
 \item We only considered a supply relationship between single pair of a source and a load. 
       Then, we set $m_\mathrm{c}=2$, which was enough under this setup. 
       At all the routers, a capacitor indexed by 1 (2) was used for forwarding (receiving) a power packet. 
 \item At the beginning of simulation, the voltage of all the capacitors within each router was set at the same value. 
       That is, at $k=1$, $v_{(i,1)} = v_{(i,2)}$ for all $i$ satisfying $a_i = \mathrm{r}$.
       Then we connected all the capacitors of a router parallel at the beginning of each time slot. 
\end{itemize}
The first and second ones reduces the factors related to the optimization. 
They are not a restriction of the proposed algorithm. 
Of course, the proposed algorithm can deal with two or more supply pairs if $m_\mathrm{c}$ increases, and with distributed values of parameters $C$, $R$ and $U$. 
The third one is to avoid the aforementioned arbitrariness of the capacitor selection algorithm.

We defined the dispatching network as $m_\mathrm{r}=9$ routers placed in a lattice structure\cite{Nawata.etal-2018}. 
Then, numerical simulation was conducted in the following three cases with different location of source(s) and a different initial voltage distribution. 

\subsection{Case I}
In this case, we assume that the load specifies the source node. 
Thus, the algorithm is going to determine the best path between the source and the load. 

\begin{figure}[tb]
 \centering
 \includegraphics[width=.5\columnwidth]{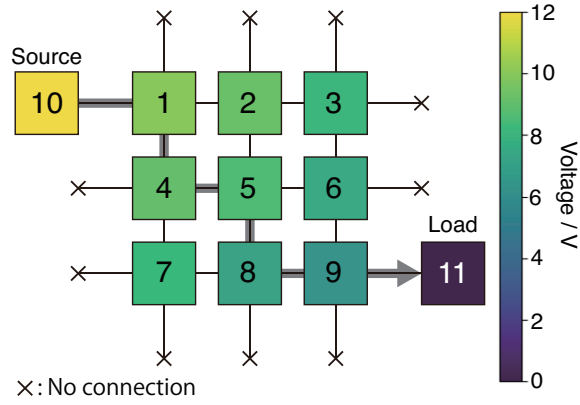}
 \caption{Result of case I. }
 \label{fig:case1}
\end{figure}

Figure~\ref{fig:case1} shows the network structure and the initial voltage distribution for case I. 
Each square and number on it indicate a node and its index, respectively. 
The color of the squares represents the voltage of each node. 
We considered the path from the source of $12.0\,\mathrm{V}$ (node 10) to the load of $10\,\Omega$ (node 11). 

As a result, the algorithm selected the path ((10,1),(1,4),(4,5),(5,8),(8,9),(9,11)) for the first power packet. 
The arrow in Fig.~\ref{fig:case1} represents the path. 
To evaluate this result, we also see the cost of other possible paths. 
\begin{figure}[tb]
 \centering
 \includegraphics[width=.5\columnwidth]{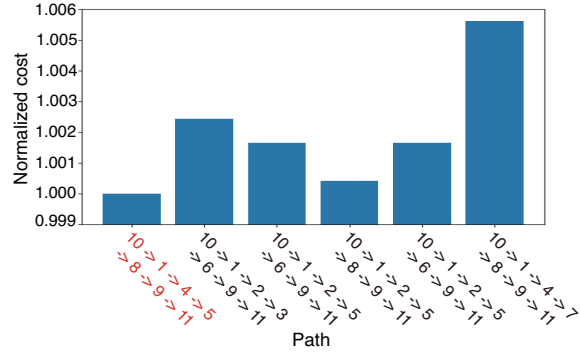}
 \caption{Cost of all the possible paths from node 10 to node 11. }
 \label{fig:cost1}
\end{figure}
Figure~\ref{fig:cost1} shows the normalized values of cost calculated along all the possible paths from node 10 to node 11. 
The values are normalized by the cost of the optimal one. 
This result indicates that the proposed algorithm successfully selected the optimal path with the lowest cost. 

Now note that each repetition of our algorithm, corresponding to the transfer of each sequential power packet, is completely independent. 
Thus, showing the result for the first power packet is adequate. 
For the subsequent packets, the same procedure is repeated. 

\subsection{Case II}
Next, we assume that the load does not specify the source node. 
The load accepts a supply from any source that can afford the requested energy $U$.

\begin{figure}[tb]
 \centering
 \includegraphics[width=.5\columnwidth]{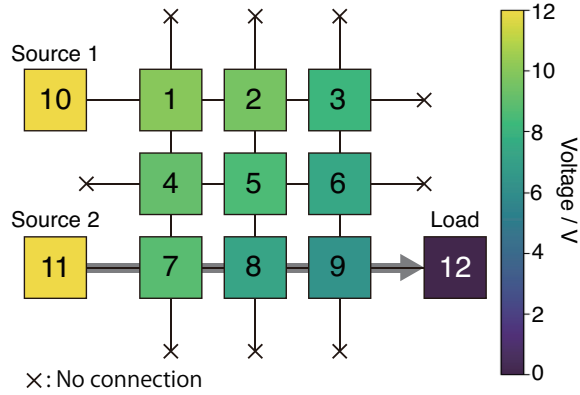}
 \caption{Result of case II. $(V_{10},V_{11})=(12,12)$. }
 \label{fig:case2}
\end{figure}

Figure~\ref{fig:case2} shows the network structure and the initial voltage distribution for case I. 
We consider the path from one of the sources of $12.0\,\mathrm{V}$ (node 10 or 11) to the load of $10\,\Omega$ (node 12). 

The algorithm selected the node 11 for the source. 
The selected path is ((11,7),(7,8),(8,9),(9,12)). 
In fact, if we fix the starting node 10, the path of the lowest cost is ((10,1),(1,4),(4,7),(7,8),(8,9),(9,12)). 
Its cost is 1.03 times as large as that of the optimal path. 
This result indicates that the proposed algorithm successfully selected the optimal source with the lower cost. 

\subsection{Case III}
The setup of this case is basically same as the case II, but we add a small disturbance to the voltage of sources 1 and 2. 
The disturbance models the fluctuation of renewable power sources such as photovoltaic cells. 

\begin{figure}[tb]
 \centering
 \includegraphics[width=.5\columnwidth]{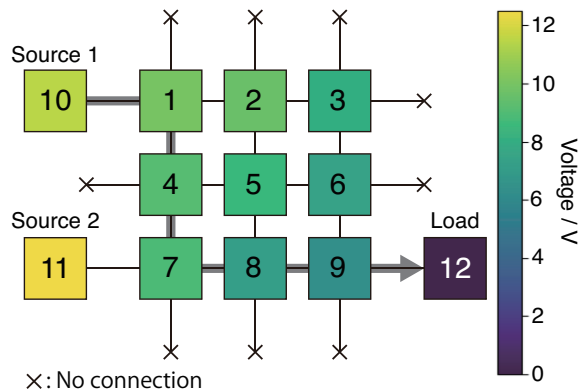}
 \caption{Result of case III. $(V_{10},V_{11})=(11.5,12.5)$.}
 \label{fig:case3}
\end{figure}

Figure~\ref{fig:case3} shows the network structure and the initial voltage distribution for case I. 
The voltage of node 10 and node 11 is set at $11.5\,\mathrm{V}$ and $12.5\,\mathrm{V}$, respectively. 
Then, similarly with case II, we consider the path from one of the sources to the load. 

The algorithm selected the node 10 for the source. 
The optimal path is ((10,1),(1,4),(4,7),(7,8),(8,9),(9,12)). 
Although there is some paths that have shorter geometric distance, e.g. ((11,7),(7,8),(8,9),(9,12)), the selected one has the lowest cost in the viewpoint of power loss. 
Since the cost varies with the ratio of node voltage, the optimality of the power flow does not necessarily correspond to that measured in geometric distance on the graph. 
This is an important aspect for considering a power flow in a dispatching system with multiple sources and loads.

\subsection{Summary of Results}
The result of case I clearly shows that the proposed scheme can find the optimal path between a specific pair of a source and a load. 
Then the results of cases II and III show that the proposed scheme can allocate the optimal source among the possible candidates to the load request. 
Of course, the path between the allocated source and the load is also optimized in the same way as in case I. 
In addition, the comparison of the results of cases II and III indicates that the optimal source differs depending on the voltage distribution of the network. 

\section{Conclusion}
In this paper, we first proposed the circuit models for connection of network nodes in a power packet dispatching system.  
Using the model, we derived the cost of power packet delivery as a function of voltage distribution. 
Then we considered the optimal routing problem based on the models and developed the algorithm to solve it. 
Through the numerical simulations we confirmed that the proposed algorithm can find the optimal path for the power delivery. 

The optimality is not measured in the sense of geometric distance but of energy loss along the path. 
Thus, the optimal path varies when the voltage distribution of the sources and the router nodes. 
As is mentioned in Introduction, the inclusion of EVs and renewable sources causes dynamic change of the voltage distribution. 
Even under this condition, the proposed method creates efficient allocation of distributed power sources to the demand of loads. 


\section*{Acknowledgments}
This work was partly supported by JSPS KAKENHI Grant Numbers 20K14732 and 20H02151, Cross-ministerial Strategic Innovation Promotion Program (SIP) 18088028, and by Program on Open Innovation Platform with Enterprises, Research Institute and Academia (OPERA) JPMJOP1841 from Japan Science and Technology Agency.

\end{document}